\documentclass[aap,MSNbibl,seceqn,dvips]{arximspdf}

\doi{10.1214/13-AAP925} 
\volume{24}
\issue{1}
\pubyear{2014}
\firstpage{312}
\lastpage{336}

\makeatletter
\newcommand{\rrVert}{\Vert}
\newcommand{\rrvert}{\vert}
\newcommand{\llVert}{\Vert}
\newcommand{\llvert}{\vert}
\newtheorem{Theorem}{Theorem}[section]
\newtheorem{Proposition}{Proposition}[section]
\newtheorem{Lemma}{Lemma}[section]
\newproclaim{Remark}{Remark}[section]
\makeatother

\begin{document}
\begin{frontmatter}

\title{A stochastic control approach to no-arbitrage bounds given marginals,
with an application to~lookback options}
\runtitle{Derivatives bounds given marginals}

\begin{aug}
\author[A]{\fnms{A.} \snm{Galichon}\ead[label=e1]{alfred.galichon@sciences-po.fr}\thanksref{t1}},
\author[B]{\fnms{P.} \snm{Henry-Labord\`ere}\ead[label=e2]{pierre.henry-labordere@sgcib.com}}
\and
\author[C]{\fnms{N.} \snm{Touzi}\corref{}\ead[label=e3]{nizar.touzi@polytechnique.edu}\thanksref{t1}}
\runauthor{A. Galichon, P. Henry-Labord\`ere and N. Touzi}
\affiliation{Sciences Po, Soci\'et\'e G\'en\'erale and Ecole Polytechnique}
\address[A]{A. Galichon\\
Department of Economics\\
Sciences Po Paris\\
27, rue Saint-Guillaume\\
75007 Paris\\
France\\
\printead{e1}}
\address[B]{P. Henry-Labord\`ere\\
Global Market Quantitative Research\\
Soci\'et\'e G\'en\'erale\\
17 cours Valmy---Paris La D\'efense 7\\
France\\
\printead{e2}}
\address[C]{N. Touzi\\
Centre de Math\'ematiques Appliqu\'ees\\
Ecole Polytechnique\\
Palaiseau, 91128\\
France\\
\printead{e3}}
\end{aug}
\thankstext{t1}{Supported by the Chair \textit{Financial Risks}
of the \textit{Risk Foundation} sponsored by Soci\'et\'e G\'en\'erale,
the Chair \textit{Derivatives of the Future} sponsored by the
F\'ed\'eration Bancaire Fran\c{c}aise, the Chair \textit{Finance
and Sustainable Development} sponsored by EDF and CA-CIB, and FiME,
Laboratoire de Finance des March\'es de l'Energie.}

\received{\smonth{9} \syear{2011}}
\revised{\smonth{2} \syear{2013}}

%
\begin{abstract}
We consider the problem of superhedging under volatility uncertainty
for an investor allowed to dynamically trade the underlying asset, and
statically trade European call options for all possible strikes with
some given maturity. This problem is classically approached by means of
the Skorohod Embedding Problem (SEP). Instead, we provide a dual
formulation which converts the superhedging problem into a continuous
martingale optimal transportation problem. We then show that this
formulation allows us to recover previously known results about
lookback options. In particular, our methodology induces a new proof of
the optimality of Az\'ema--Yor solution of the SEP for a certain class
of lookback options. Unlike the SEP technique, our approach applies to
a large class of exotics and is suitable for numerical approximation
techniques.
\end{abstract}

%

\end{frontmatter}

\section{Introduction}

In a financial market allowing for the dynamic trading of a riskless
asset, normalized to unity, and some given underlying assets without
restrictions, the fundamental theorem of asset pricing essentially
states that the absence of arbitrage opportunities is equivalent to the
existence of a probability measure under which the underlying asset
process is a martingale. See Kreps \cite{Kreps}, Harrison and Pliska
\cite{HarrisonPliska} and Delbaen and Schachermayer
\cite{DelbaenSchachermayer}. Then, for the purpose of hedging, the only
relevant information is the quadratic variation of the assets price
process under such a martingale measure. Without any further assumption
on the quadratic variation, the robust superhedging cost reduces to an
obvious bound which can be realized by static trading on the underlying
assets; see Cvitani\'c, Pham and Touzi \cite{cpt} and Frey \cite{frey}.

In this paper, we examine the problem of superhedging, under the
condition of no-arbitrage, when the financial market allows for the
static trading of European call options in addition to the dynamic
trading of the underlying asset with price process $\{X_t,t\geq0\}$.
For simplicity, we consider the case where all available European call
options have the same maturity $T$. However, we idealize the financial
market assuming that such European call options are available for all
possible strikes. Under the linearity (and continuity) assumption, this
means that although the joint distribution $\mathbb{P}$ of the process
$X$ is unknown, the modeler has access to the function $K\in\mathbb
{R}_+\longmapsto\mathbb{E}[(X_T-K)^+]$, and therefore to the marginal
distribution $\mu$ of the random variable $X_T$, as observed by Breeden
and Litzenberger \cite{bre}. In particular, any $T$-maturity vanilla
derivative, with payoff $g(X_T)$, can be perfectly replicated by a
portfolio of European calls, and has an unambiguous no-arbitrage price
$\int g\,d\mu$, as long as $g\in\mathbb{L}^1(\mu)$, which can be
expressed as a linear combination of the given prices of the underlying
calls.

This problem is classically approached in the literature by means of
the Skorohod embedding problem (SEP) which shows up naturally due to
the Dambis--Dubins--Schwarz time change result. However, only some
special cases of derivatives are eligible for this approach, namely
those defined by a payoff which is invariant under rime change. The use
of the SEP techniques to solve the robust superhedging problem can be
traced back to the seminal paper by Hobson \cite{hobson98}. The survey
paper by Hobson \cite{hobson} is very informative and contains the
relevant references on the subject. For a derivative security
$g(X_s,s\leq t)$ written on an underlying asset $X$, the idea is to
search for a function $\lambda$ and a martingale $M$ such that
$g(X_s,s\leq t)\leq\lambda(X_t)+M_t$ so that $\mathbb{E}[g(X_s,s\leq
\tau)]\leq\int\lambda\,d\mu$ for all stopping times $\tau$ such that
$X_\tau$ has the distribution $\mu$ determined from European call
options, as outlined above. The bound is then obtained by designing
$\lambda$, $M$ and $\tau$ such that
$g(X_s,s\leq\tau)=\lambda(X_\tau)+M_\tau$.

In this paper, we develop an alternative approach which relates the
robust superhedging problem to the literature on stochastic control;
see Fleming and Soner~\cite{flemingsoner}. Our approach consists of
solving directly the robust superhedging problem whose solution
provides the above function $\lambda$ and martingale $M$. As
a~consequence, unlike the SEP approach which requires that the payoff
of the exotic to be invariant under time change, there is no
restriction on the class of derivative securities as the time change
step is avoided. Moreover, our methodology is related to the optimal
transportation theory, and in fact opens the door to an original new
ramification in this theory by imposing naturally that the
transportation be performed along a continuous martingale.

Our first main result, reported in Sections~\ref{sect-superhedge} and
\ref{sect-superhedgecall}, provides a formulation of the robust
superhedging problem based on the Kantorovich duality in the spirit of
Benamou and Brenier \cite{benamoubrenier}; see Villani \cite{villani}.
A finite dimensional version of this dual formulation was already used
by Davis, Obl\'oj and Raval \cite{davisoblojraval}. Also, this duality
is implicitly proved for special cases of derivatives in the previous
literature based on the SEP. So the importance of our duality result
lies in its generality. Moreover, similar to the SEP approach, the
solution of the dual problem provides the corresponding optimal hedging
strategy and worst case model. Finally, the dual formulation is
suitable for numerical approximation techniques as shown by Bonnans and
Tan \cite{bonnanstan} and Tan and Touzi \cite{TanTouzi}.

Our next concern is to show that this duality result is not merely a
theoretical result, but may be applied to compute things in practice.
In this paper, this is demonstrated in the context of lookback
derivatives, where the robust superhedging problem is known to be
induced by the Az\'ema--Yor solution of the SEP
\mbox{\cite{AzemaYor1,AzemaYor2}}. A semi-static hedging strategy
corresponding to this bound was produced by Hobson \cite{hobson98}. Our
second main result, reported in Section~\ref{sect-Lookback}, reproduces
this bound by means of our dual formulation. In particular, this
provides a new presentation of the fact that the Az\'ema--Yor solution
of the SEP realizes the upper bound for a certain class of lookback
options. We also recover in Section~\ref{sect-forwardLookback} the
robust superhedging cost for the forward lookback option which was also
derived in Hobson~\cite{hobson98}.

The results of
Sections~\ref{sect-Lookback}~and~\ref{sect-forwardLookback} are not
new, and are reported here in order to show that the optimal
transportation approach is suitable to recover these known results.
Beyond rediscovering the formerly derived results in the literature, we
would like to insist on the fact that the present optimal
transportation approach complements the SEP approach by putting the
emphasis on the superhedging problem whose solution provides the
optimal semi-static superhedging strategy. This feature is perfectly
illustrated in the related paper \cite{HOST}, where:
\begin{longlist}[(ii)]
\item[(i)] this approach is used in order to provide the optimal
    semi-static superhedging strategy;

\item[(ii)] by using the optimal semi-static superhedging strategy from
    (i), Obl\'oj and Spoida \cite{OblojSpoida} prove an extension of
    the Az\'ema--Yor solution of the SEP to the case of multiple
    intermediate marginals.
\end{longlist}

The final outcome from \cite{HOST,OblojSpoida} is an extension of the
previous results of Brown, Hobson and Rogers \cite{brownhobsonrogers}
and Madan and Yor \cite{madanyor}.

We also observe that one should not expect to solve explicitly the
problem of robust superhedging, in general. Therefore, an important
advantage of the optimal transportation approach is that the
corresponding dual formulation is suitable for numerical approximation
techniques; see Bonnans and Tan \cite{bonnanstan} for the case of
variance options, and Tan and Touzi \cite{TanTouzi} for a class of
optimal transportation problems motivated by this paper.

Finally, we would like to emphasize that by expressing the robust bound
as the value function of a superhedging problem, we are implicitly
addressing the important issue of no-arbitrage under model uncertainty,
as discussed in the previous literature; see, for example, Cox and
Obl\'oj \cite{coxobloj} and Davis, Obl\`oj and Raval
\cite{davisoblojraval}. Indeed, if an arbitrage does exist, under some
conveniently defined weak form, the value function would be infinite.
Conversely, whenever the value function is finite, the financial market
would either admit a genuine model-free arbitrage, if existence holds
for the dual problem, or some convenient notion of weak arbitrage
otherwise.

\section{Model-free bounds of derivatives securities}\label{sect-formulation}

\subsection{The probabilistic framework}\label{sect-setup}

Let $\Omega:= \{\omega\in C([0,T], \mathbb{R}^d)\dvtx \omega_0=0\}$ be
the canonical space equipped with the uniform norm $\llVert
\omega\rrVert _\infty:= \sup_{0\leq t\leq T}|\omega_t|$, $B$~the
canonical process, $\mathbb {P}_0$ the Wiener measure, $\mathbb{F}:=
\{\mathcal{F}_t\}_{0\leq t\leq T}$ the filtration generated by $B$ and
$\mathbb{F}^+:= \{\mathcal{F}^+_t, 0\leq t\leq T\} $ the right limit
of $\mathbb{F}$, where $\mathcal{F}^+_t:= \bigcap_{s>t}
\mathcal{F}_s$.

Throughout the paper, $X_0$ is some given initial value in $\mathbb{R}_+^d$,
and we denote
\[
X_t:= X_0+B_t\qquad\mbox{for } t\in[0,T].
\]
For any $\mathbb{F}$-progressively measurable process $\alpha$ with
values in $\mathcal{S}_d^{>0}$ (space of definite positive symmetric
matrices) and satisfying $\int_0^T|\alpha_s|\,ds<\infty$,
$\mathbb{P}_0$-a.s., we define the probability measure on
$(\Omega,\mathcal{F})$,
\[
\mathbb{P}^\alpha:=\mathbb{P}_0\circ\bigl(X^\alpha
\bigr)^{-1},
\]
where
\[
X^\alpha_t:=X_0+\int_0^t \alpha_r^{1/2}\,dB_r, \qquad t
\in[0,T], \mathbb {P}_0\mbox{-a.s.}
\]
Then $X$ is a $\mathbb{P}^\alpha$-local martingale. Following
\cite{stz-aggregation}, we denote by $\overline{\mathcal{P}}_S$ the
collection of all such probability measures on $(\Omega,\mathcal{F})$.
The quadratic variation process $\langle X\rangle=\langle B\rangle$ is
universally defined under any $\mathbb{P}\in\overline{\mathcal{P}}_S$,
and takes values in the set of all nondecreasing continuous functions
from $\mathbb{R}_+$ to $\mathcal{S}_d^{>0}$. Finally, we recall from
\cite{stz-aggregation} that
%
\begin{eqnarray}\label{0-1MRP}
&\mbox{every }\mathbb{P}\in\overline{\mathcal{P}}_S \mbox{ satisfies the Blumenthal zero--one law}&
\nonumber\\[-10pt]\\[-10pt]
&\mbox{and the martingale representation property.}& \nonumber
\end{eqnarray}
In this paper, we shall focus on the subset $\mathcal{P}_{ \infty}$ of
$\overline{\mathcal{P}}_S$ consisting of all measures $\mathbb{P}$ such
that
\[
X \mbox{ is a } \mathbb{P}\mbox{-uniformly integrable martingale with values in }\mathbb{R}^d.
\]
The restriction of the probability measures in $\mathcal{P}_{ \infty}$
to those induced by uniformly integrable martingales $X$ is motivated
by our subsequent interpretation of the entries $X^i$ as price
processes of financial securities.

\subsection{Model-free super-hedging problem}

For all $\mathbb{P}\in\mathcal{P}_{ \infty}$, we denote
\[
\mathbb{H}^2_{\mathrm{loc}}(\mathbb{P}):= \biggl\{H\in
\mathbb{H}^0(\mathbb{P})\dvtx \int_0^T \operatorname{Tr}\bigl[H_t^{\mathrm{T}}H_t\,d\langle B
\rangle_t\bigr] <\infty, \mathbb{P}\mbox{-a.s.} \biggr\},
\]
where $\operatorname{Tr}$ denotes the trace operator. We assume that
\[
\mbox{\textit{the interest rate is zero.}}
\]
Under the self-financing condition, for any portfolio process $H$,
the portfolio value process
%
\begin{eqnarray}\label{XH}
Y_t^H &:=& Y_0+\int_0^t H_s\cdot dB_s,\qquad t \in[0,T],
\end{eqnarray}
is well defined $\mathbb{P}$-a.s. for every $\mathbb{P}\in\mathcal
{P}_{ \infty}$, whenever $H\in\mathbb{H}^2_{\mathrm{loc}}$. This
stochastic integral should be rather denoted ${Y^H_t}^\mathbb{P}$ to
emphasize its dependence on $\mathbb{P}$; see, however, Nutz~\cite{Nutz}.

Let $\xi$ be an $\mathcal{F}_T$-measurable random variable. We
introduce the subset of martingale measures
\[
\mathcal{P}_{ \infty}(\xi):= \bigl\{\mathbb{P}\in\mathcal{P}_{ \infty}\dvtx
\mathbb{E}^\mathbb{P}\bigl[\xi ^-\bigr]<\infty\bigr\}.
\]
The reason for restricting to this class of models is that, under the
condition that $\mathbb{E}^\mathbb{P}[\xi^+]<\infty$, the hedging
cost of $\xi$
under $\mathbb{P}$ is expected to be $-\infty$ whenever $\mathbb
{E}^\mathbb{P}[\xi
^-]=\infty$. As usual, in order to avoid doubling strategies, we
introduce the set of admissible portfolios,
\[
\mathcal{H}(\xi):= \bigl\{H\dvtx H\in\mathbb{H}^2_{\mathrm{loc}}
\mbox{ and } Y^H\mbox{ is a }\mathbb{P}\mbox{-supermartingale for all }
\mathbb {P}\in\mathcal{P}_{ \infty} (\xi) \bigr\}.
\]
The model-free superhedging problem is defined by
%
\begin{eqnarray}
\label{V0} U^0(\xi) &:=& \inf \bigl\{ Y_0\dvtx \exists H\in
\mathcal{H}(\xi), Y^H_1\geq\xi,\mathbb{P}\mbox{-a.s. for
all }\mathbb{P}\in\mathcal {P}_{ \infty} (\xi) \bigr\}.
\end{eqnarray}
We call $U^0$ the model-free superhedging bound, and we recall its
interpretation as the no-arbitrage upper bound on the market price
of the derivative security $\xi$, for an investor who has access to
continuous-time trading the underlying securities with price process
$X$.

\subsection{Dual formulation of the super-hedging bound}\label{sect-superhedge}

We denote by $\mathrm{UC}(\Omega_{X_0})$ the collection of all
uniformly continuous maps from $\Omega_{X_0}$ to $\mathbb{R}$, where
$X_0\in \mathbb{R}^d$ is a fixed initial value, and
$\Omega_{X_0}:=\{\omega\in C([0,T],\mathbb{R}
^d_+)\dvtx\omega_0=X_0\}$. The following result is a direct adaptation
from Soner, Touzi and Zhang \cite{stz-duality}.

%
\begin{Theorem}\label{thmstz}
Let $\xi\in\mathrm{UC}(\Omega_{X_0})$ be such that $\sup_{\mathbb{P}\in
\mathcal{P}_{ \infty}}\mathbb{E}^\mathbb{P}[\xi^+]<\infty$. Then
\[
U^0(\xi) = \sup_{\mathbb{P}\in\mathcal{P}_{ \infty}}\mathbb{E}^\mathbb
{P}[\xi].
\]
Assume further that $U^0(\xi)$ is finite. Then there exists a process
$H\in\mathcal{H}(\xi)$ and a family of nondecreasing predictable
processes $\{K^\mathbb{P},\mathbb{P}\in\mathcal{P}_{ \infty}(\xi)\}$,
with $K^\mathbb{P}_0=0$ for all $\mathbb{P}\in\mathcal{P}_{ \infty
}(\xi)$, s.t.
%
\begin{eqnarray}
\label{representation} \xi &=& U^0(\xi)+\int_0^1
H_t\cdot dB_t-K^\mathbb{P}_1,\qquad
\mathbb{P}\mbox{-a.s. for all }\mathbb{P}\in\mathcal{P}_{\infty}(\xi).
\end{eqnarray}
\end{Theorem}

The proof is reported in Section~\ref{sect-proof-thmduality}.

%
\begin{Remark}\label{rem-thmstz}
A similar dual representation as in Theorem \ref{thmstz} was first
obtained by Denis and Martini \cite{DM} in the bounded volatility case.
Notice, however, that the family of nondominated singular measures in
\cite{DM} is not included in our set $\overline{\mathcal{P}}_S$, and
does not allow for the existence of an optimal super-hedging strategy.

While revising this paper, a new approach for the robust superhedging
problem, initiated by Nutz and von Handel \cite{NutzvanHandel} and
Neufeld and Nutz \cite{NeufeldNutz}, allowed for an extension of
Theorem \ref{thmstz} by Possama\"{i}, Royer and Touzi
\cite{PossamaiRoyerTouzi}.
\end{Remark}

\subsection{Calibration adjusted no-arbitrage bound}\label{sect-superhedgecall}

In this section we specialize the discussion to the one-dimensional
case. This is consistent with the one-dimensional practical treatment
of vanilla options on real financial markets.

We assume that, in addition to the continuous-time trading of the
primitive securities, the investor can take static positions on
$T$-maturity European call or put options with all possible strikes
$K\geq0$. Then, from Breeden and Litzenberger \cite{bre}, the investor
can identify that the $T$-marginal distribution of the underlying asset
under the pricing measure is some probability measure $\mu\in
M(\mathbb{R})$, the set of all probability measures on $\mathbb{R}$.

%
\begin{Remark}
For the purpose of the present financial application, the measure $\mu$
has a support in $\mathbb{R}_+$. We consider, however, the general case
$\mu\in M(\mathbb{R})$ in order to compare our results to the
Az\'ema--Yor solution of the SEP.
\end{Remark}

For any scalar function $\lambda\in\mathbb{L}^1(\mu)$, the $T$-maturity
European derivative defined by the payoff $\lambda(X_T)$ has an
un-ambiguous no-arbitrage price
\[
\mu(\lambda) = \int\lambda\,d\mu
\]
and can be perfectly replicated by buying and holding a portfolio of
European calls and puts of all strikes, with the density $\lambda
''(K)$ at strike $K$ (with $\lambda''$ understood in the sense of
distributions). See Carr and Chou \cite{carrchou}. In particular, given
the spot price $X_0>0$ of the underlying assets, the probability
measure $\mu$ must satisfy
\[
\int x\mu(dx)=X_0.
\]
We now define an improvement of the no-arbitrage upper-bound by
accounting for the additional possibility of statically trading the
European call options. Let
%
\begin{eqnarray}\label{Lambdamu}
\Lambda^\mu &:=& \Bigl\{\lambda\in
\mathbb{L}^1(\mu)\dvtx \sup_{\mathbb{P}\in\mathcal{P}_{ \infty}}\mathbb{E}^\mathbb
{P}\bigl[\lambda (X_T)^-\bigr]<\infty \Bigr\}\quad\mbox{and}
\nonumber\\[-8pt]\\[-8pt]
\Lambda_{\mathrm{UC}}^\mu &:=&
\Lambda^\mu\cap\mathrm{UC}( \mathbb{R}),\nonumber
\end{eqnarray}
where $\mathrm{UC}(\mathbb{R})$ is the collection of all uniformly
continuous maps from $\mathbb{R}$ to $\mathbb{R}$. For all
$\lambda\in\Lambda^\mu$, we denote $\xi ^\lambda:=\xi-\lambda(X_1)$.
The improved no-arbitrage upper bound is defined by
%
%
\begin{eqnarray}
\label{V} U^\mu(\xi) &:=& \inf \bigl\{ Y_0\dvtx \exists \lambda
\in\Lambda_{\mathrm{UC}}^\mu \mbox{ and }H\in\mathcal{H}\bigl(
\xi^\lambda\bigr),
\nonumber\\[-8pt]\\[-8pt]
&&\hspace*{18pt}\overline{Y}{}^{H,\lambda}_1\geq\xi,\mathbb{P}\mbox{-a.s. for all }\mathbb{P}\in\mathcal{P}_{ \infty
}\bigl(\xi^\lambda\bigr) \bigr\}, \nonumber
\end{eqnarray}
where $\overline{Y}{}^{H,\lambda}$ denotes the portfolio value of a
self-financing strategy with continuous trading $H$ in the primitive
securities, and static trading $\lambda$ in the $T$-maturity European
calls with all strikes
%
\begin{eqnarray}
\label{barX} \overline{Y}{}^{H,\lambda}_1 &:=& Y^H_1-
\mu(\lambda)+\lambda(X_T),
\end{eqnarray}
indicating that the investor has the possibility of buying at time
$0$ any derivative security with payoff $\lambda(X_T)$ for the price
$\mu(\lambda)$.

The next result is a direct application of Theorem \ref{thmstz}.

%
\begin{Proposition}\label{propVximudual}
\hspace*{-1.5pt}Let $\mu\in M(\mathbb{R})$, and
$\xi\in\mathrm{UC}(\Omega_{X_0})$ with\break
$\sup_{\mathbb{P}\in\mathcal{P}_{ \infty}}\mathbb{E}^\mathbb
{P}[\xi^+]<\infty$. Then
\[
U^\mu(\xi) = \inf_{\lambda\in\Lambda_{\mathrm{UC}}^\mu} \sup_{\mathbb{P}\in\mathcal{P}_{ \infty}}
\bigl\{\mu(\lambda)+\mathbb{E}^\mathbb{P} \bigl[\xi-\lambda
(X_T) \bigr] \bigr\}.
\]
\end{Proposition}

\begin{pf}
Observe that
\[
U^\mu(\xi) = \inf_{\lambda\in\Lambda_{\mathrm{UC}}^\mu} U^0 \bigl(
\xi+\mu(\lambda)-\lambda(X_T) \bigr).
\]
For every fixed $\lambda$, if $V(0):=\sup_{\mathbb{P}\in\mathcal{P}_{
\infty}}\mathbb{E}^\mathbb{P}[\xi+\mu(\lambda)-\lambda (X_T)]<\infty$,
then the proof of Theorem \ref{thmstz}, reported in
Section~\ref{sect-proof-thmduality}, applies and we get $U^0
(\xi+\mu(\lambda)-\lambda(X_T) )=V(0)$. On the other hand, if
$V(0)=\infty$, then notice from the proof of Theorem \ref{thmstz} that
the inequality $U^0 (\xi+\mu(\lambda)-\lambda(X_T) )\geq V(0)$ is still
valid in this case, and therefore $U^0 (\xi+\mu(\lambda)-\lambda(X_T)
)= V(0)$.
\end{pf}

%
\begin{Remark} \label{rem-sanitycheck}
As a sanity check, let us consider the case $\xi=g(X_T)$,\break  for some
uniformly continuous function $g$ with $\mu(|g|)<\infty$ and\break
$\sup_{\mathbb{P}\in\mathcal{P}_{ \infty} }\mathbb{E}^\mathbb{P}
[|g(X_T)|]<\infty$, and let us verify that $U^\mu(\xi)=\mu(g)$.

First, since $g\in\Lambda_{\mathrm{UC}}^\mu$, we may take $\lambda=g$,
and it follows from the dual formulation of Proposition
\ref{propVximudual} that $U^\mu(\xi)\leq\mu(g)$. On the other hand, it
is easily seen that
\[
\sup_{\mathbb{P}\in\mathcal{P}_{ \infty}}\mathbb{E}^\mathbb{P}\bigl[g(X_T)
\bigr] = g^{\mathrm{conc}}(X_0),
\]
where $g^{\mathrm{conc}}$ is the smallest concave majorant of $g$.
Then, it follows from the dual formulation of Proposition
\ref{propVximudual} that
$U^\mu(\xi)=\inf_{\lambda\in\Lambda_{\mathrm{UC}}^\mu}\mu
(\lambda)+(g-\lambda)^{\mathrm{conc}}(X_0)\geq\inf_{\lambda\in\Lambda_{\mathrm{UC}}^\mu}\mu
(\lambda)+\mu(g-\lambda)=\mu(g)$ as expected.
\end{Remark}

%
\begin{Remark}\label{rem-solution}
Similar to the SEP approach, the dual formulation of Proposition
\ref{propVximudual} gives access to the optimal hedging strategy and
the worst case model. This requires that we prove an additional
existence result of a solution to the inf-sup problem $(\lambda
^*,\mathbb{P} ^*)$. Then $\lambda^*$ is the optimal $T$-maturity
vanilla profile, and $\mathbb{P}^*$ is the worst case model
corresponding to the upper bound. The optimal dynamic hedging strategy
in the underlying asset is, as usual, obtained by representation of the
residual security $\xi-\lambda^*(X_T)$; see the proof of Theorem
\ref{thmstz} in Section~\ref{sect-proof-thmduality}.
\end{Remark}

%
\begin{Remark}\label{numerics}
The dual formulation of Proposition \ref{propVximudual} is suitable for
numerical approximation. Indeed, for each fixed multiplier $\lambda$,
the maximization problem is a (singular) stochastic control problem
which may be approximated by finite differences or Monte Carlo methods.
Then optimization stage with respect to $\lambda$ requires an
additional iteration. This issue is addressed in Tan and Bonnans
\cite{bonnanstan} and Tan and Touzi \cite{TanTouzi}.
\end{Remark}

\subsection{Connection with optimal transportation theory}

As an alternative point of view, one may directly imbed in the
no-arbitrage bounds the calibration constraint that the risk neutral
marginal distribution of $B_T$ is given by $\mu$.

For convenience of comparison with the optimal transportation theory,
the discussion of this subsection will be focused on the no-arbitrage
lower bound. A~natural formulation of the calibration adjusted
no-arbitrage lower bound is
%
\begin{eqnarray}
\label{Uoptimaltransport}
\ell(\xi,\mu) &:=& \inf \bigl\{\mathbb{E}^\mathbb{P}[
\xi]\dvtx\mathbb{P}\in\mathcal {P}_{ \infty},X_0
\sim_\mathbb{P}\delta_{X_0} \mbox{ and } X_T
\sim_\mathbb{P}\mu \bigr\},
\end{eqnarray}
where $\delta_{X_0}$ denotes the Dirac mass at the point $X_0$. We
observe that a direct proof that $\ell(\xi,\mu)$ coincides with the
corresponding sub-hedging cost is not obvious in the present context.

Under this form, the problem appears as minimizing the coupling
criterion $E^\mathbb{P}[\xi]$ which involves the law of the process $X$
under $\mathbb{P}$, over all those probability measures $\mathbb
{P}\in\mathcal{P}_{ \infty}$ such that the marginal distributions of
$X$ at times $0$~and~$T$ are fixed. This is the general scope of
optimal transportation problems as introduced by Monge and Kantorovich;
see, for example, Villani \cite{villani} and Mikami and Thieulen
\cite{mikamithieulen}. Motivated by the present financial application,
Tan and Touzi \cite{TanTouzi} extended the Kantorovich duality as
described below. However, the above problem $\ell(\xi,\mu)$ does not
satisfy the assumptions in \cite{TanTouzi} so that none of the results
contained in this literature apply to our context.

The classical approach in optimal transportation consists of deriving a
dual formulation for problem (\ref{Uoptimaltransport}) by means of the
classical convex duality theory. Recall that $M(\mathbb{R}_+)$ denotes
the collection of all probability measures on $\mathbb{R}_+$. Then, the
Legendre dual with respect to $\mu$ is defined by
\[
\ell^*(\xi,\lambda):= \sup_{\mu\in M(\mathbb{R}_+)} \bigl\{\mu(\lambda)-\ell(
\xi,\mu ) \bigr\}\qquad\mbox{for all } \lambda\in C^0_b(\mathbb{R}_+)
\]
the set of all bounded continuous functions from $\mathbb{R}_+$ to
$\mathbb{R}$.
Direct calculation shows that
\begin{eqnarray*}
\ell^*(\xi,\lambda) &=& \sup \bigl\{\mathbb{E}^\mathbb{P}\bigl[
\lambda(X_T)-\xi\bigr]\dvtx \mu\in M(\mathbb{R}_+), \mathbb{P}\in
\mathcal{P}_{ \infty},X_0\sim_\mathbb{P} \delta_{X_0} \mbox{ and }
X_T\sim_\mathbb{P}\mu \bigr \}
\\
&=& \sup \bigl\{\mathbb{E}^\mathbb{P}\bigl[\lambda(X_T)-\xi
\bigr]\dvtx \mathbb{P}\in\mathcal{P}_{ \infty},X_0
\sim_\mathbb{P}\delta_{X_0} \bigr\}.
\end{eqnarray*}
Observe that the latter problem is a standard (singular) diffusion
control problem.

It is easily checked that $\ell$ is convex in $\mu$. However, due to
the absence of a uniform bound on the quadratic variation of $X$ under
$\mathbb{P}\in\mathcal{P}_{ \infty}$, it is not obvious whether it is
lower semicontinuous with respect to $\mu$. If the latter property were
true, then the equality $\ell^{**}=\ell$ provides
\[
\ell(\xi,\mu) = \sup_{\lambda\in C^0_b} \bigl\{\mu(\lambda)-\ell^*(\xi,
\lambda) \bigr\},
\]
which is formally (up to the spaces choices) the lower bound analogue
of the dual formulation of Proposition \ref{propVximudual}. A
discrete-time analysis of this duality is contained in the parallel
work to the present one by Beiglb\"ock, Henry-Labord\`ere and
Penkner~\cite{beiglbockhenrylaborderepenkner}. We also observe that,
for special cases of payoffs $\xi$, this duality was implicitly proved
in the previous literature based on the SEP approach; see, for example,
Cox, Hobson and Obl\'oj \cite{coxhobsonobloj}.

\section{Application to lookback derivatives}\label{sect-Lookback}

Throughout this section, we consider the one dimensional case $d=1$.
The derivative security is defined by the lookback payoff
%
\begin{eqnarray}
\label{xilookback} \xi &=& g \bigl(X^*_T \bigr)\qquad\mbox{where }
X^*_T:=\max_{t\leq T} X_t
\end{eqnarray}
and
%
\begin{eqnarray}
\label{hyp-g}
g\dvtx\mathbb{R}&\longrightarrow&\mathbb{R}_+ \mbox{ is a
$C^1$ nondecreasing function.}
\end{eqnarray}
Our main interest is to show that the optimal upper bound given by
Proposition~\ref{propVximudual},
\[
U^\mu(\xi) = \inf_{\lambda\in\Lambda_{\mathrm{UC}}^\mu} \bigl\{ \mu(
\lambda)+u^\lambda(0,X_0,X_0) \bigr\}
\]
reproduces the already known bound corresponding to the Az\'ema--Yor
solution to the Skorohod embedding problem. Here, $u^\lambda$ is the
value function of the dynamic version of stochastic control problem
%
\begin{eqnarray}
\label{stochcontrollookback}
u^\lambda(t,x,m) &:=& \sup_{\mathbb{P}\in\mathcal{P}_{ \infty}}
\mathbb{E}^\mathbb{P} \bigl[g \bigl(M^{t,x,m}_T\bigr)-
\lambda\bigl(X^{t,x}_T\bigr) \bigr],\qquad t\leq T, (x,m)\in\bolds{\Delta},\hspace*{-15pt}
\end{eqnarray}
where $\bolds{\Delta}:= \{(x,m)\in\mathbb{R}^2\dvtx x\leq m \}$, and
\[
X_u^{t,x}:=x+(B_u-B_t),\qquad
M_u^{t,x,m}:=m\vee\max_{t\leq r\leq u}X_r^{t,x},\qquad
0\leq t\leq u\leq T.
\]
When the time origin is zero, we shall simply write $X^x_u:=X^{0,x}_u$
and $M^{x,m}_u:=M^{0,x,m}_u$.

For the subsequent analysis, we also claim that
%
\begin{eqnarray}
\label{U-max} U^\mu(\xi) &=& \inf_{\lambda\in\Lambda^\mu} \bigl\{ \mu(
\lambda)+u^\lambda(0,X_0,X_0) \bigr\},
\end{eqnarray}
where $\Lambda^\mu$ is defined in (\ref{Lambdamu}). This follows from
the recent extension by \cite{PossamaiRoyerTouzi}, which appeared
during the revision of this paper (see Remark \ref{rem-thmstz}), and
avoids placing further conditions on $\mu$ to ensure that the function
$\lambda^*$ defined in (\ref{lambda*}) below is uniformly continuous.

\subsection{Formulation in terms of optimal stopping}
We first convert the optimization problem $u^\lambda$ into an infinite
horizon optimal stopping problem.

%
\begin{Proposition}\label{prop-max-optimalstop}
For any $\lambda\in\Lambda^\mu$, the functions $u^\lambda$ is
independent of $t$ and
%
\begin{eqnarray}
\label{rep-max-optimalstop} u^\lambda(x,m) &=& \sup_{\tau\in\mathcal{T}_{ \infty}}
\mathbb{E}^{\mathbb{P}_0} \bigl[g\bigl(M^{x,m}_{\tau}\bigr)-
\lambda\bigl(X^x_{\tau}\bigr) \bigr]\qquad\mbox{for all } (x,m)\in\bolds{\Delta},
\end{eqnarray}
where $\mathcal{T}_{ \infty}$ is the collection of all stopping times
$\tau$ such that the stopped process $\{X_{t\wedge\tau},t\geq0\}$ is a
$\mathbb{P} _0$-uniformly integrable martingale.
\end{Proposition}

\begin{pf}
The present argument is classical, but we
could not find a clear
reference. We therefore report it for completeness.\vadjust{\goodbreak}

By the definition of $\mathcal{P}_{ \infty}$, we may write the
stochastic control problem (\ref{stochcontrollookback}) in its strong
formulation
\[
u^\lambda(t,x,m):= \sup_{\sigma\in\Sigma^+} \mathbb{E}^{\mathbb{P}_0}
\bigl[g \bigl(M^{\sigma,t,x,m}_T \bigr)-\lambda\bigl(X^{\sigma,t,x}_T
\bigr) \mid \bigl(X^{\sigma,t,x}_t,M^{\sigma,t,x,m}_t
\bigr)=(x,m) \bigr],
\]
where
\[
X_s^{\sigma,t,x}=x+\int_t^s
\sigma_r\,dB_r,\qquad M_s^{\sigma,t,x,m}:=m\vee
\max_{t\leq r\leq s} X^{\sigma,t,x}_r,\qquad 0\leq t\leq s\leq T
\]
and $\Sigma^+$ is the set of all nonnegative progressively measurable
processes, with $\int_0^T\sigma_s^2\,ds<\infty$, and such that the
process $\{X_s^{\sigma,t,x},t\leq s\leq T\}$ is a uniformly integrable
martingale.

We shall denote $\phi(x,m):=g(m)-\lambda(x)$.
\begin{longlist}[(2)]
\item[(1)] For a stopping time $\tau\in\mathcal{T}_\infty$, we
    define the processes
    $\sigma^\tau_t:=\mathbf{1}_{\{\tau\wedge(t/(T-t))\}}$ and
    $X^{\sigma^\tau}_t:=\int_0^t\sigma^\tau_s\,dB_s$, $t\in[t,T]$. Then
    the corresponding measure $\mathbb{P}^{\sigma^\tau}:=\mathbb
    {P}_0\circ(X^{\sigma^\tau})^{-1}\in\mathcal{P}_{ \infty}$, and
%
\begin{eqnarray}
\label{DDSineq1} \sup_{\tau\in\mathcal{T}_\infty}\mathbb{E}^{\mathbb{P}_0}\bigl[\phi
(X_\tau,M_\tau)\bigr] &\leq& \sup_{\mathbb{P}\in\mathcal{P}_{ \infty}}
\mathbb{E}^{\mathbb
{P}}\bigl[\phi(X_\tau,M_\tau)\bigr].
\end{eqnarray}

\item[(2)] To obtain the reverse inequality, we observe by the
    Dambis--Dubins--Schwarz theorem (see, e.g., Karatzas and Shreve
    \cite{KaratzasShreve}, Theorem 3.4.6) that the law of $(X,M)$ under
    $\mathbb{P}=\mathbb{P} ^\sigma\in\mathcal{P}_{ \infty}$ is the same
    as the law of $(X_\tau,M_\tau)$ under $\mathbb{P}_0$ where
    $\tau:=\int_0^T\sigma^2_t\,dt$ is a stopping time with respected to
    the time-changed filtration $\{\mathcal{F}^B_{T_t},t\geq0\} $,
    $T_t:=\inf\{s\dvtx\langle B\rangle_s >t\}$. In order to convert to
    the context of the canonical filtration of $B$, we use the result
    of Szpirglas and Mazziotto \cite{szpirglas-mazziotto}. This allows
    us to conclude that equality holds in (\ref{DDSineq1}).\qed
\end{longlist}\noqed
\end{pf}

In view of the previous results, we are reduced to the problem
%
\begin{eqnarray}
\label{boundlookback} U^\mu(\xi) &:=& \inf_{\lambda\in\Lambda_0^\mu}
\bigl\{\mu(\lambda)+u^\lambda(X_0,X_0) \bigr\},
\end{eqnarray}
where
\[
u^\lambda(X_0,X_0):= \sup_{\tau\in\mathcal{T}_{ \infty}}
J(\lambda,\tau),\qquad J(\lambda,\tau):= \mathbb{E}^{\mathbb{P}_0} \bigl[g
\bigl(X^*_\tau\bigr)-\lambda(X_\tau) \bigr]
\]
and the set $\Lambda^\mu$ of (\ref{Lambdamu}) translates in the present
context to
%
\begin{eqnarray}
\label{Lambda0mu} \Lambda_0^\mu &=& \Bigl\{\lambda\in
\mathbb{L}^1(\mu)\dvtx \sup_{\tau\in\mathcal{T}_{ \infty}}\mathbb{E}\bigl[
\lambda(X_\tau )^-\bigr]<\infty \Bigr\}.
\end{eqnarray}

\subsection{The main result}

The endpoints of the support of the distribution $\mu$ are denoted by
\[
\ell^\mu:=\sup \bigl\{x\dvtx\mu \bigl([x,\infty) \bigr)=1 \bigr\}
\quad\mbox{and}\quad
r^\mu:=\inf \bigl\{x\dvtx\mu \bigl((x,\infty) \bigr)=0 \bigr\}.\vadjust{\goodbreak}
\]
The Az\'ema--Yor solution of the \textit{Skorohod embedding problem} is
defined by means of the so-called barycenter function,
%
\begin{eqnarray}
\label{barycenter} b(x) &:=& \frac{\int_{[x,\infty)} y\mu(dy)} {
\mu ([x,\infty) )} \mathbf{1}_{\{x<r^\mu\}} +x \mathbf{1}_{\{x\geq r^\mu\}},\qquad x\geq0.
\end{eqnarray}

%
\begin{Remark}\label{rem-barycenter}
Hobson \cite{hobson98} observed that the barycenter function can be
alternatively defined as the left-continuous inverse to the following
function $\beta$. Given the European calls prices $c(x):=\int
(y-x)^+\mu(dy)$ and $X_0=\int y\mu(dy)$, define the function
%
\begin{eqnarray}
\beta(x)&:=&\max \biggl\{\mathop{\arg\min}_{y<x}\frac{c(y)}{x-y} \biggr\}\qquad\mbox{for } x\in[X_0,r^\mu),
\\
\beta(x)&=&\ell^\mu\qquad\mbox{for } x\in[0,X_0)\quad\mbox{and}\quad\beta(x)=x\qquad\mbox{for } x\in[r^\mu,\infty).\hspace*{-20pt}
\end{eqnarray}
On $[X_0,r^\mu)$, $\beta(x)$ is the largest minimizer of the function
$y\longmapsto c(y)/(x-y)$ on $(-\infty,x)$. Then, $\beta$ is
nondecreasing, right-continuous, and $\beta(x)<x$ for all $x\in
[X_0,r^\mu)$. Notice that $\beta(X_0)=\ell^\mu:=\sup\{x\dvtx\mu ((0,x]
)>0\}$.
\end{Remark}

Finally, we introduce the Hardy--Littlewood transform
$\mu^{\mathrm{HL}}$ of $\mu$,
%
\begin{eqnarray}
\label{muHL}
\mu^{\mathrm{HL}} \bigl([y,\infty) \bigr) &:=& \inf
_{\xi<y} \frac{c(y)}{y-\xi} = \frac{c (\beta(y) )}{y-\beta(y)}\qquad\mbox{for all } y\geq0,
\end{eqnarray}
where the functions $c$ and $\beta$ are defined in the previous remark;
see Proposition~4.10(c) in Carraro, El Karoui and Obl\'oj
\cite{carraro-elkaroui-obloj}.

The following result is a combination of \cite{hobson98,peskir-max}.
Our objective is to derive it directly from the dual formulation of
Proposition \ref{propVximudual}. Let
%
\begin{eqnarray}
\label{tau*} \tau^* &:=& \inf \bigl\{t>0\dvtx X^*_t\geq b(X_t) \bigr\}
\end{eqnarray}
and
%
\begin{eqnarray}
\label{lambda*} \lambda^*(x) &:=& \int_{\ell^\mu}^x\int
_{\ell^\mu}^y g' \bigl(b(\xi) \bigr)
\frac{b(d\xi)}{b(\xi)-\xi}\,dy;\qquad x<r^\mu.
\end{eqnarray}
Notice that $\lambda^*\in[0,\infty]$ as the integral of a nonnegative
function. To see that \mbox{$\lambda^*<\infty$}, we compute by the
Fubini theorem that
\[
\lambda^*(x) = \int_{\ell^\mu}^x g'
\bigl(b(\xi) \bigr) \frac{x-\xi}{b(\xi)-\xi} b(d\xi)
\]
and we observe that $(x-\xi)/(b(\xi)-\xi)$ is bounded near $\ell^\mu$.
Then $\lambda^*(x)\leq C(x)[g(b(x))-g(X_0)]<\infty$ for some constant
$C(x)$ depending on $x$.

%
\begin{Theorem}\label{thmlookback}
Let $\mu\in M(\mathbb{R})$, and $\xi=g(X^*_T)$ for some $C^1$
nondecreasing function $g$ satisfying $\sup_{\mathbb{P}\in\mathcal{P}_{
\infty}}\mathbb{E}^\mathbb{P}[\xi^+]<\infty$, and
$\mu^{\mathrm{HL}}(g)<\infty$. Then
\[
U^\mu(\xi) = \mu\bigl(\lambda^*\bigr)+J\bigl(\lambda^*,\tau^*\bigr)= \mu^{\mathrm{HL}}(g).
\]
\end{Theorem}

The proof is reported in the subsequent section.

\subsection{An upper bound for the optimal upper bound}

In this section, we prove that
%
\begin{eqnarray}
\label{upperbound} U^\mu(\xi) &\leq& \mu\bigl(\lambda^*\bigr)+J\bigl(\lambda^*,\tau^*\bigr).
\end{eqnarray}
Our first step is to use the following construction due to Peskir
\cite{peskir-max} which provides a guess of the value function
$u^\lambda$ for functions $\lambda$ in the subset
%
\begin{eqnarray}
\hat\Lambda_0^\mu &:=& \bigl\{\lambda\in
\Lambda_0^\mu\dvtx\lambda\mbox{ is convex}\bigr\}.
\end{eqnarray}
By classical tools from stochastic control theory, the value function
$u^\lambda(x,m)$ is expected to solve the dynamic programming equation
%
\begin{eqnarray}\label{DPE}
\min \bigl\{u^\lambda-g+\lambda,-u^\lambda_{xx}\bigr\}&=&0\qquad\mbox{on }
\bolds{\Delta}\quad\mbox{and}
\nonumber\\[-8pt]\\[-8pt]
u^\lambda_m(m,m)&=&0\qquad\mbox{for }m\in\mathbb{R}.\nonumber
\end{eqnarray}
The first part of the above DPE is an ODE for which $m$ appears only as
a parameter involved in the domain on which the ODE must hold. Since we
are restricting to convex $\lambda$, one can guess a solution of the form
%
\begin{eqnarray}
\label{vpsi0} v^\psi(x,m) &:=& g(m)-\lambda \bigl(x\wedge\psi(m) \bigr) -
\lambda' \bigl(\psi(m) \bigr) \bigl(x-x\wedge\psi(m) \bigr),
\end{eqnarray}
that is, $v^\psi(x,m)=g(m)-\lambda(x)$ for $x\leq\psi(m)$ and is given
by the tangent at the point $\psi(m)$ for $x\in[\psi(m),m]$. For later
use, we observe that for $x\in[\psi(m),m]$,
%
\begin{eqnarray}\label{vpsi}
v^\psi(x,m) &=& g(m)-\lambda\bigl(\psi(m)\bigr) +\int
_{\psi(m)}^x \frac{\partial}{\partial y}\bigl\{
\lambda'(y) (x-y)\bigr\}\,dy
\nonumber\\[-8pt]\\[-8pt]
&=& g(m)-\lambda(x) +\int_{\psi(m)}^x (x-y)
\lambda''(dy)\qquad\mbox{for } x\in\bigl[\psi(m),m\bigr],\hspace*{-30pt}\nonumber
\end{eqnarray}
where $\lambda''$ is the second derivative measure of the convex
function $\lambda$.

We next choose the function $\psi$ in order to satisfy the Neumann
condition in~(\ref{DPE}). Assuming that $\lambda$ is smooth, we obtain
by direct calculation that the free boundary $\psi$ must verify the
ordinary differential equation (ODE)
%
\begin{eqnarray}
\label{ODE0} \lambda'' \bigl(\psi(m) \bigr)
\psi'(m) &=& \frac{g'(m)}{m-\psi(m)}\qquad\mbox{for all } m\in\mathbb{R}.
\end{eqnarray}
For technical reasons, we need to consider this ODE in the relaxed
sense. This contrasts our analysis with that of Peskir
\cite{peskir-max} and Obl\'oj \cite{Obloj-max}. Since $\lambda$ is
convex, its second derivative $\lambda''$ is well defined as measure on
$\mathbb{R}_+$. We then introduce the weak formulation of the ODE
(\ref{ODE0}),
%
\begin{eqnarray}
\label{ODE} \int_{\psi(B)}\lambda''(dy)
&=& \int_B\frac{g'(m)}{m-\psi(m)}\,dm\qquad\mbox{for all } B\in\mathcal{B}(\mathbb{R})
\end{eqnarray}
and we introduce the collection of all relaxed solutions of
(\ref{ODE0}),
%
\begin{eqnarray}\label{Psilambda}
\Psi^\lambda &:=& \bigl\{\psi\mbox{ right-continuous:
(\ref{ODE}) holds and }
\nonumber\\[-8pt]\\[-8pt]
&&\hspace*{59pt}\psi(m)<m\mbox{ for all }m\in\mathbb{R} \bigr\}.\nonumber
\end{eqnarray}

%
\begin{Remark}\label{rem-ODEaffine}
For later use, we observe that (\ref{ODE}) implies that all functions
$\psi\in\Psi^\lambda$ are nondecreasing. Indeed, for $y_1\leq y_2$, it
follows from (\ref{ODE}), together with the nondecrease of $g$ in
(\ref{hyp-g}) and the convexity of $\lambda$, that
\begin{eqnarray*}
\psi(y_2) &=& \bigl(\lambda'\bigr)^{-1}
\biggl(\lambda' \bigl(\psi(y_1)+ \bigr) +\int
_{y_1}^{y_2}\frac{g'(m)} {
m-\psi(m)}\,dm \biggr)
\\
&\geq& \bigl(
\lambda'\bigr)^{-1} \bigl(\lambda' \bigl(
\psi(y_1)+ \bigr) \bigr) \geq \psi(y_1),
\end{eqnarray*}
where $(\lambda')^{-1}$ is the right-continuous inverse of the
nondecreasing function $\lambda'$.
Then, by direct integration that
\[
\mbox{the function } x\longmapsto \lambda(x)-\int_{X_0}^x
\int_{X_0}^{\psi^{-1}(y)}\frac{g'(\xi)}{\xi-\psi(\xi)}\,d\xi\,dy \mbox{ is
affine,}
\]
where $\psi^{-1}$ is the right-continuous inverse of $\psi$. This
follows from direct differentiation of the above function in the sense
of generalized derivatives.
\end{Remark}

A remarkable feature of the present problem is that there is no natural
boundary condition for the ODE (\ref{ODE0}) or its relaxation
(\ref{ODE}). The following result extends the easy part of the elegant
maximality principle proved in Peskir \cite{peskir-max} by allowing for
possibly nonsmooth functions $\lambda$. We emphasize the fact that our
approach does not need the full strength of Peskir's maximality
principle.

%
\begin{Lemma}\label{lemPeskir}
Let $\lambda\in\hat{\Lambda}_0^\mu$ and $\psi\in\Psi^\lambda$ be
arbitrary. Then $u^\lambda\leq v^\psi$.
\end{Lemma}

\begin{pf}
We organize the proof in three steps:

\begin{longlist}[(2)]
\item[(1)] We first prove that $v^\psi$ is differentiable
    in $m$ on the diagonal with
%
\begin{eqnarray}
\label{Neumann} v^\psi_m(m,m)&=&0\qquad\mbox{for all }m\in \mathbb{R}.
\end{eqnarray}
Indeed, since $\psi\in\Psi_\lambda$, it follows from Remark
\ref{rem-ODEaffine} that
\[
\lambda(x) = c_0+c_1x +\int_{X_0}^x
\int_{X_0}^{\psi^{-1}(y)}\frac{g'(\xi)}{\xi-\psi(\xi)}\,d\xi\,dy
\]
for some scalar constants $c_0$, $c_1$. Plugging this expression into
(\ref{vpsi0}), we see that
\begin{eqnarray*}
v^\psi(x,m) &=& g(m)- \biggl(c_0+c_1\psi(m) +
\int_{X_0}^{\psi(m)}\int_{X_0}^{\psi^{-1}(y)}
\frac{g'(\xi)}{\xi-\psi(\xi)}\,d\xi\,dy \biggr)
\\
&&{} - \biggl(c_1+\int_{X_0}^m\frac{g'(\xi)}{\xi-\psi(\xi)}\,d\xi \biggr) \bigl(x-\psi(m) \bigr)
\\
&=& g(m)-c_0-c_1x+\int_{X_0}^m
\frac{g'(\xi)}{\xi-\psi(\xi)} \bigl(\psi(\xi)-x\bigr)\,d\xi,
\end{eqnarray*}
where the last equality follows from the Fubini theorem together with
the fact that $g$ is nondecreasing and $\psi(\xi)<\xi$. Since $g$ is
differentiable, (\ref{Neumann}) follows by direct differentiation with
respect to $m$.

\item[(2)] For an arbitrary stopping time
    $\tau\in\mathcal{T}_{ \infty}$, we introduce the stopping times
    $\tau_n:=\tau\wedge\inf\{ t>0\dvtx\llvert X_t-x\rrvert >n\}$. Since
    $v^\psi$ is concave in $x$, as a consequence of the convexity of
    $\lambda$, it follows from the It\^o--Tanaka formula that
\begin{eqnarray*}
v^\psi(x,m) &\geq& v^\psi(X_{\tau_n},M_{\tau_n})
-\int_0^{\tau_n} v^\psi_x(X_t,M_t)\,dB_t -\int_0^{\tau_n}
v^\psi_{m}(X_t,M_t)\,dM_t
\\
&\geq& g(M_{\tau_n})-\lambda(X_{\tau_n}) -\int
_0^{\tau_n} v^\psi_x(X_t,M_t)\,dB_t
-\int_0^{\tau_n} v^\psi_{m}(X_t,M_t)\,dM_t
\end{eqnarray*}
by the fact that $v^\psi\geq g-\lambda$. Notice that
$(M_t-X_t)\,dM_t=0$. Then by the Neumann condition (\ref{Neumann}), we
have $v^\psi_{m}(X_t,M_t)\,dM_t=v^\psi_{m}(M_t,M_t)\,dM_t=0$. Taking
expectations in the last inequality, we see that
%
\begin{eqnarray}
\label{ineq-avantlimit} v^\psi(x,m) &\geq& \mathbb{E}_{x,m}
\bigl[g(M_{\tau_n})-\lambda(X_{\tau_n}) \bigr].
\end{eqnarray}

\item[(3)] We finally take the limit as $n\to\infty$ in
    the
    last inequality. First, recall that $(X_{t\wedge\tau})_{t\geq0}$ is
    a uniformly integrable martingale. Then, by the Jensen inequality,
    $\lambda(X_{\tau_n})\leq\mathbb{E}[\lambda(X_\tau)\mid \mathcal
    {F}_{\tau_n}]$. Since
    $\lambda(X_\tau)^-\in\mathbb{L}^1(\mathbb{P}^0)$, this implies that
    $\mathbb{E}[\lambda (X_{\tau_n})]\leq\mathbb{E}[\lambda(X_\tau)]$
    where we also used the tower property of conditional expectations.
    We then deduce from (\ref{ineq-avantlimit}) that
\[
v^\psi(x,m) \geq \lim_{n\to\infty}\mathbb{E}_{x,m}
\bigl[g(M_{\tau_n})-\lambda (X_\tau ) \bigr] =
\mathbb{E}_{x,m} \bigl[g(M_{\tau})-\lambda(X_\tau)
\bigr]
\]
\end{longlist}
by the nondecrease of the process $M$ and the function $g$ together
with the monotone convergence theorem. By the arbitrariness of $\tau
\in\mathcal{T}_{ \infty}$, the last inequality shows that $v^\psi\geq
u^\lambda$.
\end{pf}

Our next result involves the function
%
\begin{eqnarray}\label{varphi}
\varphi(x,m):=\frac{c(x)-c_0(x)\mathbf{1}_{m<X_0}} {m-x}
\nonumber\\[-8pt]\\[-8pt]
\eqntext{\mbox{with }c_0(x):=(X_0-x)^+,(x,m)\in\bolds{\Delta}}
\end{eqnarray}
and we recall that $c(x):=\int(\xi-x)^+\mu(d\xi)$ is the (given)
European call price with strike $x$.

%
\begin{Lemma}\label{lemmu+ulambda}
For $\lambda\in\hat{\Lambda}_0^\mu$ and $\psi\in\Psi^\lambda$,
we have
\[
\mu(\lambda)+u^\lambda(X_0,X_0) \leq
g(X_0)+\int\varphi \bigl(\psi(m),m \bigr)g'(m)\,dm.\vadjust{\goodbreak}
\]
\end{Lemma}

\begin{pf}
(1) Let $\alpha\in\mathbb{R}_+$ be an arbitrary point of
differentiability of $\lambda$. Then
\[
\lambda(x) = \lambda(\alpha)+\lambda'(\alpha) (x-\alpha) +\int
_\alpha^x (x-y)\lambda''(dy).
\]
Integrating with respect to $\mu-\delta_{X_0}$ and taking
$\alpha<X_0$, this provides
\begin{eqnarray*}
&& \mu(\lambda)-\lambda(X_0)
\\
&&\qquad  = \lambda'(\alpha) \biggl( \int
x\mu(dx)-X_0 \biggr)+\int \biggl(\int_\alpha^x (x-y)\lambda''(dy)
\biggr) (\mu-\delta_{X_0}) (dx)
\\
&&\qquad = -\int_\alpha^{X_0}(X_0-y)
\lambda''(dy) +\int\mathbf{1}_{\{x\geq\alpha\}}\int
_\alpha^x (x-y)^+\lambda''(dy)\mu(dx)
\\
&&\quad\qquad{} +\int\mathbf{1}_{\{x<\alpha\}}\int_x^\alpha(y-x)
\lambda''(dy) \mu(dx).
\end{eqnarray*}
Then sending $\alpha$ to $\ell^\mu$, it follows from the convexity of
$\lambda$ together with the monotone convergence theorem that
\[
\mu(\lambda)-\lambda(X_0) = \int(c-c_0) (y)
\lambda''(dy).
\]

(2) By the inequality in Lemma \ref{lemPeskir}, together with
(\ref{vpsi}), we now compute that
\begin{eqnarray*}
&& \mu(\lambda)+u^\lambda(X_0,X_0)
\\
&&\qquad \leq
g(X_0) +\int \bigl(c(y)-c_0(y) (1_{\{y<X_0\}}
-1_{\{\psi(X_0)<y<X_0\}}) \bigr)\lambda''(dy)
\\
&&\qquad = g(X_0) +\int \bigl(c(y)-c_0(y)1_{\{y<\psi(X_0)\}}
\bigr)\lambda''(dy).
\end{eqnarray*}
We next use the ODE (\ref{ODE0}) satisfied by $\psi$ in the
distribution sense. This provides
\[
\mu(\lambda)+u^\lambda(X_0,X_0) \leq
g(X_0) +\int\frac{c(\psi(m))-c_0(\psi(m))\mathbf{1}_{\{m<X_0\}}} {
m-\psi(m)} g'(m)\,dm.
\]
Here, we observe that the endpoints in the last integral can be taken
to $0$ and $\infty$ by the nonnegativity of the integrand.
\end{pf}

We now have all ingredients to express the upper bound
(\ref{upperbound}) explicitly in terms of the barycenter function $b$
of (\ref{barycenter}).

%
\begin{Lemma}\label{lem-upperbound}
For a nondecreasing $C^1$ payoff function $g$, we have
\[
\inf_{\lambda\in\Lambda_0^\mu} \bigl\{\mu(\lambda)+u^\lambda(X_0,X_0)
\bigr\} \leq \mu^{\mathrm{HL}}(g).
\]
\end{Lemma}

\begin{pf}
Since $\hat{\Lambda}_0^\mu\subset
\Lambda_0^\mu$, we compute from Lemma \ref{lemmu+ulambda} that
%
\begin{eqnarray}\label{upperboundStep0}
&& \inf_{\lambda\in\Lambda_0^\mu} \bigl\{\mu(\lambda)+u^\lambda(X_0,X_0)
\bigr\}\nonumber
\\
&&\qquad \leq \inf_{\lambda\in\hat{\Lambda}_0^\mu} \bigl\{\mu(\lambda)+u^\lambda(X_0,X_0)\bigr\}
\\
&&\qquad \leq g(X_0)+ \inf_{\lambda\in\hat{\Lambda}_0^\mu} \inf
_{\psi\in\Psi^\lambda} \int\varphi \bigl(\psi(m),m \bigr)g'(m)\,dm.\nonumber
\end{eqnarray}
In the next two steps, we prove that the last minimization problem on
the right-hand side of (\ref{upperboundStep0}) can be solved by
pointwise minimization inside the integral. Then, in step~(3), we
compute the induced upper bound.

\begin{longlist}[(3)]
\item[(1)] For all $\lambda\in\hat{\Lambda}_0^\mu$ and
    $\psi \in\Psi^\lambda$,
\[
\int\varphi \bigl(\psi(m),m \bigr)g'(m)\,dm \geq \int\inf
_{\xi<m}\varphi (\xi,m )g'(m)\,dm.
\]
Observe that $c(x)\geq c_0(x)$ for all $x\geq0$, and
$\lim_{x\to0}c(x)-c_0(x)=0$. Then
%
\begin{eqnarray}
\label{linearissubopt1} \inf_{\xi<m} \varphi(\xi,m)&=&
\varphi(0,m)=0\qquad\mbox{for } m<X_0.
\end{eqnarray}
On the other hand, it follows from Remark \ref{rem-barycenter} that
%
\begin{eqnarray}
\label{linearissubopt2} \inf_{\xi<m} \varphi(\xi,m)&=&\inf
_{\xi<m} \frac{c(\xi)}{m-\xi} = \frac{c (\beta(m) )}{m-\beta(m)}\qquad\mbox{for } m \geq X_0.
\end{eqnarray}
By (\ref{linearissubopt1}) and (\ref{linearissubopt2}), we obtain the
lower bound
\[
\int\varphi \bigl(\psi(m),m \bigr)g'(m)\,dm \geq \int\varphi
\bigl(\beta(m),m \bigr)g'(m)\,dm.
\]

\item[(2)] We now observe that the function $\beta$, obtained by
    pointwise minimization in the previous step, solves the ODE
    (\ref{ODE}). Therefore, in order to complete the proof, it remains
    to verify that $\lambda^*\in\hat{\Lambda}_0^\mu$. The convexity of
    $\lambda^*$ is obvious. Also, since $\lambda^*\geq0$, we only need
    to prove that $\lambda^*\in\mathbb{L}^1(\mu)$. By step~(1) of the
    proof of Lemma \ref{lemmu+ulambda}, we are reduced to verifying
    that $\int c(x)(\lambda ^*)''(dx)<\infty$. Since, by definition,
    $\lambda^*$ satisfies the ODE (\ref{ODE}) with $\psi=b^{-1}$, we
    directly compute that
\begin{eqnarray*}
\int c(x) \bigl(\lambda^*\bigr)''(dx) &=& \int
\frac{c (b^{-1}(m) )} {
m-b^{-1}(m)} g'(m)\,dm
\\
&=& \int g'(m)
\mu^{\mathrm{HL}} \bigl([m,\infty) \bigr)\,dm <\infty
\end{eqnarray*}
by our assumption that $\mu^{\mathrm{HL}}(g)<\infty$.

\item[(3)] From (\ref{upperboundStep0}) and the previous two
    steps, we have
\begin{eqnarray*}
\inf_{\lambda\in\Lambda_0^\mu}
\bigl\{\mu(\lambda)+u^\lambda(X_0,X_0) \bigr\}
&\leq& g(X_0)+\int\frac{c(\beta(x))}{x-\beta(x)}g'(x)\,dx
\\
&=& g(X_0)+\int\mu^{\mathrm{HL}} \bigl([y,\infty) \bigr)
g'(x)\,dx = \mu^{\mathrm{HL}}(g)
\end{eqnarray*}
\end{longlist}
by a direct integration by parts.
\end{pf}

\subsection{\texorpdfstring{Completing the proof of Theorem \protect\ref{thmlookback}}
{Completing the proof of Theorem 3.1}}\label{sect-proofequality}

To complete the proof of the theorem, it remains to prove that
\[
\inf_{\lambda\in\Lambda_0^\mu} \bigl\{\mu(\lambda)+u^\lambda(X_0,X_0)
\bigr\} \geq \mu^{\mathrm{HL}}(g).
\]
To see this, we use the fact that the stopping time $\tau^*$ defined in
(\ref{tau*}) is a solution of the Skorohod embedding problem, that is,
$X_{\tau^*}\sim\mu$ and $(X_{t\wedge\tau^*})_{t\geq0}$ is a uniformly
integrable martingale; see Az\'ema and Yor \cite{AzemaYor1,AzemaYor2}.
Moreover $X^*_{\tau^*}\sim\mu^{\mathrm{HL}}$. Then, for all
$\lambda\in\Lambda_0^\mu$, it follows from the definition of
$u^\lambda$ that $u^\lambda(X_0,X_0)\geq J(\lambda,\tau^*)$, and
therefore
\begin{eqnarray*}
\mu(\lambda)+u^\lambda(X_0,X_0) &\geq& \mu(
\lambda) +\mathbb{E}_{X_0,X_0} \bigl[g\bigl(X^*_{\tau^*}\bigr)-
\lambda(X_{\tau^*}) \bigr]
\\
&=& \mathbb{E}_{X_0,X_0} \bigl[g
\bigl(X^*_{\tau^*}\bigr) \bigr] = \mu^{\mathrm{HL}}(g).
\end{eqnarray*}

\section{Forward start lookback options}\label{sect-forwardLookback}

In this section, we provide a second application to the case where the
derivative security is defined by the payoff
\[
\xi = g \bigl(B^*_{t_1,t_2} \bigr)\qquad\mbox{where } B^*_{t_1,t_2}:=\max_{t_1\leq t\leq t_2}B_t
\]
and $g$ satisfies the same conditions as in the previous section. We
assume that the prices of call options $c_1(k)$ and $c_2(k)$ for the
maturities $t_1$ and $t_2$ are given for all strikes,
\[
c_1(k)=\int(x-k)^+\mu_1(dx)\quad\mbox{and}\quad
c_2(k)=\int(x-k)^+\mu_2(dx),\qquad k\geq0.
\]
We also assume that $\mu_1\preceq\mu_2$ are in convex order:
\[
c_1(0)=c_2(0) \quad\mbox{and}\quad c_1(k)\leq
c_2(k)\qquad\mbox{for all }k\geq0.
\]
The model-free superhedging cost is defined as the minimal initial
capital which allows to superhedge the payoff $\xi$, quasi-surely, by
means of some dynamic trading strategy in the underlying stock, and a
static strategy in the calls $(c_1(k))_{k\geq0}$ and
$(c_2(k))_{k\geq0}$.

This problem was solved in Hobson \cite{hobson98} in the case $g(x)=x$.
Our objective here is to recover his results by means of our stochastic
control approach.

A direct adaptation of Proposition \ref{propVximudual} provides the
dual formulation of this problem as
\[
U^{\mu_1,\mu_2}(\xi) = \sup_{(\lambda_1,\lambda_2)
\in\Lambda^{\mu_1}\times\Lambda^{\mu_2}} \mu_1(
\lambda_1)+\mu_2(\lambda_2)
+u^{\lambda_1,\lambda_2}(X_0,X_0),
\]
where
\[
u^{\lambda_1,\lambda_2}(x,m):= \sup_{\mathbb{P}\in\mathcal{P}_{ \infty}} \mathbb{E}^\mathbb{P}_{x,m}
\bigl[g\bigl(B^*_{t_1,t_2}\bigr)-\lambda_1(B_{t_1}) -
\lambda_2(B_{t_2}) \bigr].
\]
We next observe that the dynamic value function corresponding to the
stochastic control problem $u^{\lambda_1,\lambda_2}$ reduces to our previously
studied problem $u^{\lambda_2}$ at time $t_1$. Then, it follows from
the dynamic programming principle that
\begin{eqnarray*}
U^{\mu_1,\mu_2}(\xi) &=& \inf_{(\lambda_1,\lambda_2)\in\Lambda^{\mu_1}\times\Lambda^{\mu_2}} \mu_1(
\lambda_1)+\mu_2(\lambda_2)
\\
&&{} +\sup
_{\mathbb{P}\in\mathcal{P}_{ \infty}}\mathbb{E}^\mathbb {P} \bigl[u^{\lambda
_2}(B_{t_1},B_{t_1})-\lambda_1(B_{t_1}) \bigr].
\end{eqnarray*}
Since the expression to be maximized only involves the distribution of
$B_{t_1}$, it follows from Remark \ref{rem-sanitycheck} together with
the Dambis--Dubins--Schwarz time change formula that
\[
U^{\mu_1,\mu_2}(\xi) = \inf_{\lambda_2\in\Lambda^{\mu_2}_0} \mu_2(
\lambda_2) +\int u^{\lambda_2}(x,x)\mu_1(dx).
\]
We next obtain an upper bound by restricting attention to the subset
$\hat\Lambda_0^{\mu_2}$ of convex multipliers of $\Lambda^{\mu_2}_0$.
For such multipliers, we use the inequality $u^{\lambda_2}\leq
v^{\psi_2}$ for all $\psi_2\in\Psi^{\lambda_2}$ as derived in Lemma
\ref{lemPeskir}. This provides
\begin{eqnarray*}
&& U^{\mu_1,\mu_2}(\xi)
\\
&&\qquad \leq
\inf_{\lambda_2\in\hat{\Lambda}^{\mu_2}_0} \mu_2( \lambda_2) +\int
v^{\psi_2}(x,x)\mu_1(dx)
\\
&&\qquad = \mu_1(g)+ \inf_{\lambda_2\in\hat{\Lambda}^{\mu_2}_0} \inf
_{\psi_2\in\Psi^\lambda} \mu_2(\lambda_2)-
\mu_1(\lambda_2)
 +\int\int_{\psi_2(x)}^x
(x-y)\lambda_2''(dy)\mu_1(dx)
\\
&&\qquad = \mu_1(g)+ \inf_{\lambda_2\in\hat{\Lambda}^{\mu_2}_0} \inf
_{\psi_2\in\Psi^\lambda}\int \biggl( c_2(y)-c_1(y)
\\
&&\hspace*{146pt}{}+ \int(x-y)\mathbf{1}_{\{\psi_2(x)<y<x\}}\mu_1(dx) \biggr)\lambda_2''(dy)
\\
&&\qquad = \mu_1(g)+ \inf_{\lambda_2\in\hat{\Lambda}^{\mu_2}_0} \inf
_{\psi_2\in\Psi^\lambda} \int \biggl( c_2(y) -\int(x-y)\mathbf{1}_{\{y\leq\psi_2(x)\}}\mu_1(dx) \biggr)\lambda_2''(dy)
\\
&&\qquad = \mu_1(g)
+ \inf_{\lambda_2\in\hat{\Lambda}^{\mu_2}_0} \inf
_{\psi_2\in\Psi^\lambda}
\int \biggl( c_2\bigl(
\psi_2(m)\bigr) -\int\bigl(x-\psi_2(m)\bigr)\mathbf{1}_{\{m\leq x\}}\mu_1(dx) \biggr)
\\[-2pt]
&&\hspace*{139pt}
{}\times \frac{g'(m)\,dm}{m-\psi_2(m)}
\\[-2pt]
&&\qquad = \mu_1(g)+ \inf_{\lambda_2\in\hat{\Lambda}^{\mu_2}_0} \inf
_{\psi_2\in\Psi^\lambda} \int \biggl(\frac{c_2 (\psi_2(m) )-c_1(m)} {
m-\psi_2(m)} -\mu_1 \bigl([m,\infty) \bigr) \biggr)g'(m)\,dm,
\end{eqnarray*}
where the last equalities follow from similar manipulations as in Lemma
\ref{lemmu+ulambda}, and in particular make use of the ODE (\ref{ODE}).
Since $g'\geq0$, we may prove, as in the case of lookback options, that
the above minimization problem reduces to the pointwise minimization of
the integrand, so that the optimal obstacle is given by
\[
\psi_2^*(x) = \max \Bigl\{\mathop{\arg\min}_{\xi<x} h(\xi) \Bigr\}\qquad\mbox{where } h(\xi):= \frac{c_2(\xi)-c_1(m)} {
m-\xi}, \xi<m.
\]
Notice that $h$ has left and right derivative at every $\xi<m$, with
\[
h'(\xi) = \frac{c_2(\xi)+(x-\xi)c_2'(\xi)-c_1(x)} {
(x-\xi)^2},\qquad\mbox{a.e.},
\]
where the numerator is a nondecreasing function of $\xi$, takes the
positive value $c_2(x)-c_1(x)$ at $\xi=x$, and takes the negative value
$X_0-x-c_1(x)$ at $\xi=0$. Then $\psi_2^*(x)$ is the largest root of
the equation
%
\begin{eqnarray}
\label{hatpsiforwardlookback} c_2 \bigl(\psi^*_2(x)
\bigr) + \bigl(x-\psi^*_2(x) \bigr)c_2'
\bigl(\psi^*_2(x) \bigr) &=& c_1(x),\qquad\mbox{a.e.}
\end{eqnarray}
so that $h$ is nonincreasing to the left of $\psi^*_2(m)$ and
nondecreasing to its right.

At this point, we recognize exactly the solution derived by Hobson
\cite{hobson98}. In particular, $\psi^*_2$ induces a solution
$\tau_2^*$ to the Skorohod embedding problem, and we may use the
expression of $u^{\lambda_2}$ as the value function of an optimal
stopping problem. Then, we may conclude the proof that the upper bound
derived above is the optimal upper bound by arguing as in
Section~\ref{sect-proofequality} that
\[
u^{\lambda_2}(x,x) \geq \mathbb{E}_{x,x} \bigl[g
\bigl(X^*_{\tau_2^*}\bigr)-\lambda_2(X_{\tau_2^*}) \bigr].
\]
We get that the upper bound is given by
%
\begin{eqnarray}\label{boundforwardlookback}
U^{\mu_1,\mu_2}(\xi) &=& \mu_1(g)-\int g'(m)\mu_1 \bigl([m,\infty) \bigr)\,dm\nonumber
\\[-2pt]
&&{} +\int \biggl(\frac{c_2 (\psi^*_2(m) )-c_1(m)} {m-\psi^*_2(m)}\biggr)g'(x)\,dx
\nonumber\\[-9pt]\\[-9pt]
&=& \mu_1(g) -\int \bigl(c_2' \bigl(\psi^*_2(m) \bigr)-c_1'(m)\bigr)g'(m)\,dm\nonumber
\\[-2pt]
&=& g\bigl(\ell^\mu\bigr) -\int c_2' \bigl(\psi^*_2(m) \bigr)g'(m)\,dm\nonumber
\end{eqnarray}
by (\ref{hatpsiforwardlookback}).\vadjust{\goodbreak}

\section{Proof of the duality result}\label{sect-proof-thmduality}

Let $\xi\dvtx\Omega\longrightarrow\mathbb{R}$ be a measurable map with
$\sup_{\mathbb{P}\in\mathcal{P}_{ \infty}}\mathbb{E}^\mathbb
{P}[\xi^+]<\infty$. If $\mathcal{P}_{ \infty}(\xi)=\varnothing$, the
result is trivial. We then continue assuming that $\mathcal{P}_{
\infty}(\xi)\neq\varnothing$ and therefore $U_0(\xi)>-\infty$. Let
$X_0\in\mathbb{R}$ be such that
%
\begin{equation}\label{X1gexi}
X_T^H\geq\xi\qquad\mbox{for some }H\in \mathcal{H}.
\end{equation}
By definition of the admissibility set $\mathcal{H}(\xi)$, it follows
that the process $X^H$ is a $\mathbb{P}$-local martingale and a
$\mathbb{P}$-supermartingale for any $\mathbb{P}\in\mathcal{P}_{
\infty}(\xi)$. Then, it follows from (\ref{X1gexi}) that
$X_0\geq\mathbb{E}^\mathbb{P}[\xi]$ for all $\mathbb{P}\in\mathcal
{P}_{ \infty} (\xi)$. From the arbitrariness of $X_0$ and $\mathbb{P}$,
this shows that
%
\begin{equation}\label{firstinequality}
U^0(\xi) \geq
\sup_{\mathbb{P}\in\mathcal{P}_{ \infty}(\xi )} \mathbb
{E}^\mathbb{P}[\xi] = \sup_{\mathbb{P}\in\mathcal{P}_{ \infty}}\mathbb
{E}^\mathbb{P}[\xi].
\end{equation}
In the subsequent sections, we prove that the converse inequality holds
under the additional requirement that
$\xi\in\mathrm{UC}(\Omega_{X_0})$. Following \cite{stz-duality}, this
result is obtained by introducing a \textit{dynamic version} of the
problem which is then proved to have a decomposition leading to the
required result. Due to the fact that family of probability measures
$\mathcal{P}_{ \infty}$ is nondominated, we need to define conditional
distributions on all of the probability space without excepting any
zero measure set.

\subsection{Regular conditional probability distribution}

Let $\mathbb{P}$ be an arbitrary probability measure on $\Omega$, and
$\tau$ be an $\mathbb{F}$-stopping time. The regular conditional
probability distribution (r.c.p.d.) $\mathbb{P}^\omega_\tau$ is defined
by:
\begin{longlist}[--]
\item[--] for all $\omega\in\Omega$, $\mathbb{P}^\omega_\tau$ is a
    probability measure on $\mathcal{F}_T$;

\item[--] for all $E\in\mathcal{F}_T$, the mapping
    $\omega\longmapsto\mathbb{P}^\omega_\tau(E)$ is $\mathcal{F}_\tau
    $-measurable;

\item[--] for every bounded $\mathcal{F}_T$-measurable random variable
    $\xi$, we have $\mathbb{E}^\mathbb{P}[\xi\mid
    \mathcal{F}_\tau](\omega)=\mathbb
    {E}^{\mathbb{P}^\omega_\tau}[\xi]$, $\mathbb{P}$-a.s.;

\item[--] for all $\omega\in\Omega$, $\mathbb{P}^\omega_\tau
    [\omega'\in
    \Omega\dvtx \omega'(s)=\omega(s), 0\leq s\leq\tau(\omega) ]=1$.
\end{longlist}

The existence of the r.c.p.d. is justified in Stroock and Varadhan
\cite{SV}. For a better understanding of this notion, we introduce the
shifted canonical space
\[
\Omega^t:= \bigl\{\omega\in C\bigl([t,T], \mathbb{R}^d
\bigr)\dvtx \omega(t)=0\bigr\}\qquad\mbox{for all } t\in[0,T],
\]
we denote by $B^t$ the shifted canonical process on $\Omega^t$,
$\mathbb{P}^t_0$ the shifted Wiener measure and $\mathbb{F}^{t}$ the
shifted filtration generated by $B^t$. For $0\leq s\leq t\leq T$ and
$\omega\in\Omega^s$:

\begin{longlist}[--]
\item[--] the shifted path $\omega^t\in\Omega^t$ is defined by
\[
\omega^t_r:= \omega_r-
\omega_t\qquad\mbox{for all } r\in[t, T];
\]

\item[--] the concatenation path
    $\omega\otimes_{t}\tilde\omega\in\Omega^s$, for some
    $\tilde\omega\in\Omega^t$, is defined by
\[
(\omega\otimes_t\tilde\omega) (r):= \omega_r\mathbf{1}_{[s,t)}(r) + (\omega_{t} + \tilde\omega_r){\mathbf 1}_{[t, T]}(r)\qquad\mbox{for all } r\in[s,T];
\]

\item[--] the shifted $\mathcal{F}^{t}_{T}$-measurable r.v.
    $\xi^{t,\omega }$ of some $\mathcal{F}^{s}_{T}$-measurable r.v.
    $\xi$ on $\Omega^s$ is defined by
\[
\xi^{t, \omega}(\tilde\omega):=\xi(\omega\otimes_t \tilde
\omega)\qquad\mbox{for all } \tilde\omega\in\Omega^t.
\]
\end{longlist}

Similarly, for an $\mathbb{F}^{s}$-progressively measurable process $X$
on $[s, T]$, the shifted process $\{X^{t, \omega}_r, r\in[t,T]\}$ is
$\mathbb{F}^{t}$-progressively measurable.

For notational simplicity, we set
\[
\omega\otimes_\tau\tilde\omega:= \omega\otimes_{\tau(\omega)} \tilde\omega,\qquad
\xi^{\tau,\omega}:= \xi^{\tau(\omega),\omega},\qquad
X^{\tau,\omega}:=X^{\tau(\omega),\omega}.
\]
The r.c.p.d. $\mathbb{P}^\omega_\tau$ induces a probability measure
$\mathbb{P}^{\tau,\omega}$ on $\mathcal{F}^{\tau(\omega)}_T$ such that
the \mbox{$\mathbb{P}^{\tau,\omega}$-}distribution of $B^{\tau(\omega)}$ is
equal to the $\mathbb{P}^\omega_\tau$-distribution of
$\{B_t-B_{\tau(\omega)}, t\in [\tau(\omega),T]\}$. Then, the r.c.p.d.
can be understood by the identity
\[
\mathbb{E}^{\mathbb{P}^\omega_\tau}[\xi] = \mathbb{E}^{\mathbb{P}^{\tau,\omega}}\bigl[
\xi^{\tau,\omega}\bigr]\qquad\mbox{for all } \mathcal{F}_T\mbox{-measurable r.v. }\xi.
\]
We shall also call $\mathbb{P}^{\tau,\omega}$ the r.c.p.d. of
$\mathbb{P}$.

For $0\leq t\leq T$, we follow the same construction as in
Section~\ref{sect-setup} to define the martingale measures $\mathbb
{P}^{t,\alpha}$ for each $\mathbb{F}^{t}$-progressively measurable
$\mathcal{S}^{>0}_d$-valued process $\alpha$ such that
$\int_{t}^{T}\llvert \alpha_r\rrvert\,dr <\infty$,
$\mathbb{P}^{t}_0$-a.s. The collection of all such measures is denoted
$\overline\mathcal{P}{}_S^{t}$. The subset $\mathcal{P}^t_\infty$ and the
density process $\hat a^{t}$ of the quadratic variation process
$\langle B^{t}\rangle$ are also defined similarly.

\subsection{The duality result for uniformly continuous payoffs}

Since $\xi\in\break \mathrm{UC}(\Omega_{X_0})$, there exists a modulus of
continuity function $\rho$ such that for all $t\in[0,T]$ and $\omega,
\omega'\in\Omega$, $\tilde\omega\in \Omega^t$,
\[
\bigl|\xi^{t,\omega}(\tilde\omega) - \xi^{t,\omega'}(\tilde \omega) \bigr| \leq
\rho\bigl(\bigl\|\omega-\omega'\bigr\|_t\bigr),
\]
where $\|\omega\|_t:= \sup_{0\leq s\leq t}|\omega_s|$, $0\leq t\leq T$.
The main object in the present proof is the following dynamic value
process:
%
\begin{eqnarray}
\label{def-V} V_t(\omega) &:=&\sup_{\mathbb{P}\in\mathcal{P}_{ \infty}^t} \mathbb
{E}^{\mathbb{P}_t^\omega
}[\xi]\qquad\mbox{for all } (t,\omega)\in[0,T]\times\Omega.
\end{eqnarray}
It follows from the uniform continuity property of $\xi$ that
%
\begin{equation}\label{V-adapted}
\bigl\{V_t,t\in[0,T]\bigr\}\mbox{ is a right-continuous
$\mathbb{F}$-adapted process.}
\end{equation}
Moreover, since $\sup_{\mathbb{P}\in\mathcal{P}_{ \infty}}\mathbb
{E}^\mathbb{P}[\xi^+]<\infty$, it follows that for all
$\mathbb{P}\in\mathcal{P}_{ \infty}(\xi)$ that $V_t\in
\mathbb{L}^1(\mathbb{P})$, and by following exactly the proof of
Proposition 4.7 in \cite{stz-duality}, we see that $\{V_t,t\in[0,T]\}$
is a $\mathbb{P}$-supermartingale. We may then apply the Doob--Meyer
decomposition, and deduce the existence of a pair of processes
$(H^\mathbb{P},K^\mathbb{P})$, with
$H^\mathbb{P}\in\mathbb{H}^2_{\mathrm{loc}}(\mathbb{P})$ and
$K^\mathbb{P}$ $\mathbb{P}$-integrable nondecreasing, such that
\[
V_t = V_0+\int_0^t
H^\mathbb{P}_s\,dB_s-K^\mathbb{P}_t,\qquad
t \in[0,T],\mathbb {P}\mbox{-a.s.}
\]
Since $V$ is a right-continuous semimartingale under each
$\mathbb{P}\in\mathcal{P}_{ \infty}(\xi)$, it follows from Karandikar
\cite{Karandikar} that the family of processes
$\{H^\mathbb{P},\mathbb{P}\in\mathcal{P}_{ \infty}(\xi)\} $ (defined
\mbox{$\mathbb{P}$-}a.s.) can be aggregated into a process $\hat H$ defined on
$[0,T]\times\Omega$ by $d\langle V, B\rangle_t = \hat H_t\,d\langle
B\rangle_t$, in the sense that $\hat H=H^\mathbb{P}$, $dt\times
d\mathbb{P}$-a.s. for all $\mathbb{P}\in\mathcal{P}_{ \infty}(\xi)$.
Thus we have
\[
V_t = V_0+\int_0^t\hat H_s\,dB_s-K^\mathbb{P}_t,\qquad
t\in[0,T], \mathbb{P}\mbox{-a.s. for all }\mathbb{P}\in\mathcal{P}_{ \infty
}(\xi).
\]
With $X_0:=V_0$, we see that:
\begin{longlist}[--]
\item[--] the process $X^{\hat H}:=X_0+\int_0^.\hat H_s\,dB_s$ is
    bounded from below by $V$ which is in turn bounded from below by
    $M_t^\mathbb{P}:=\mathbb{E}_t^\mathbb{P}[\xi]$, $t\in[0,T]$; since
    $\xi\in\mathbb{L}^1(\mathbb{P})$, the latter is a
    $\mathbb{P}$-martingale; consequently, $X^{\hat H}$ is a
    $\mathbb{P} $-supermartingale for all $\mathbb{P}\in\mathcal{P}_{
    \infty}(\xi)$,

\item[--] and $X^{\hat
    H}_T=V_T+K^\mathbb{P}_T=\xi+K^\mathbb{P}_T\geq\xi$,
    $\mathbb{P}$-a.s. for every $\mathbb{P}\in\mathcal{P}_{
    \infty}(\xi)$.
\end{longlist}

Then $V_0\geq U^0(\xi)$ by the definition of $U^0$.

Notice that, as a consequence of the supermartingale property of
$X^{\hat H}$ under every $\mathbb{P}\in\mathcal{P}_{ \infty}(\xi)$, we
have
\[
V_0+\sup_{\mathbb{P}\in\mathcal{P}_{ \infty}(\xi)}\mathbb {E}^\mathbb{P}
\bigl[-K^\mathbb{P}_T \bigr] \geq \sup_{\mathbb{P}\in\mathcal{P}_{ \infty}(\xi)}
\mathbb {E}^\mathbb{P} \bigl[X^{\hat H}_T-K^\mathbb{P}_T
\bigr] = \sup_{\mathbb{P}\in\mathcal{P}_{ \infty}(\xi)}\mathbb {E}^\mathbb{P}[\xi] =
V_0.
\]
Since $K^\mathbb{P}_0=0$ and $K^\mathbb{P}$ is nondecreasing, this
implies that
\[
X^{\hat H} \mbox{ is a $\mathbb{P}$-martingale for all } \mathbb{P}\in
\mathcal{P}_{ \infty}(\xi)
\]
and the nondecreasing process $K^\mathbb{P}$ satisfies the minimality
condition
\[
\inf_{\mathbb{P}\in\mathcal{P}_{ \infty}(\xi)}\mathbb {E}^\mathbb{P} \bigl[K^\mathbb{P}_T
\bigr] = 0.
\]

%
\begin{Remark}\label{rem-weak1}
A possible extension of Theorem \ref{thmstz} can be obtained for a
larger class of payoff functions $\xi$. Indeed, notice that the uniform
continuity assumption on $\xi$ is essentially used to obtain the
measurability (\ref{V-adapted}) of dynamic value process $V$ defined in
(\ref{def-V}). In the context of the application of
Section~\ref{sect-Lookback}, such an extension is needed in order to
avoid restricting our framework to those measures $\mu$ which induce a
uniformly continuous optimal static hedging $\lambda^*$. A~convenient
extension of Theorem \ref{thmstz}, relaxing the uniform integrability
condition, is obtained in Possama\"{i}, Royer and Touzi
\cite{PossamaiRoyerTouzi}, building on the recent results of Nutz and
van Handel \cite{NutzvanHandel} and Neufeld and Nutz
\cite{NeufeldNutz}.
\end{Remark}


\section*{Acknowledgments}
This version has benefited from detailed comments of two anonymous referees.
We are sincerely grateful to them for their time and effort and, in particular, for pointing out some deficiencies of the
first version.



\printaddresses

\end{document}